\documentclass[reqno,12pt]{amsart}
\usepackage{amsthm}
\usepackage{amsmath,mathrsfs}
\usepackage{color}
\usepackage{amssymb}
\usepackage{hyperref}
\usepackage{enumerate}
\usepackage{latexsym,array}
\usepackage{amsfonts}
\usepackage{shadow}

\newtheorem{Pa}{Paper}[section]
\newtheorem{Tm}[Pa]{{\bf Theorem}}
\newtheorem{La}[Pa]{{\bf Lemma}}
\newtheorem{Dn}[Pa]{{\bf Definition}}
\newtheorem{Cy}[Pa]{{\bf Corollary}}

\newtheorem{Rk}[Pa]{{\bf Remark}}
\newtheorem{Pn}[Pa]{{\bf Proposition}}

\def\hh{\mathbb{H}}

\newcommand{\rr}{\mathbb{R}}

\author[D. Alpay]{Daniel Alpay}
\address{(DA) Department of Mathematics\\
Ben-Gurion University of the Negev\\
Beer-Sheva 84105 Israel} \email{dany@math.bgu.ac.il}

\author[F. Colombo]{Fabrizio Colombo}
\address{(FC) Politecnico di
Milano\\Dipartimento di Matematica\\Via E. Bonardi, 9\\20133
Milano, Italy}
\email{fabrizio.colombo@polimi.it}

\author[J. Gantner]{Jonathan Gantner}
\address{(JG)
Vienna University of Technology\\
Institute for Analysis and Scientific Computing\\
Wiedner Hauptstrasse 8 - 10\\
1040 Wien, Austria
} \email{jonathan.gantner@gmx.at}

\author[I. Sabadini]{Irene Sabadini}
\address{(IS) Politecnico di
Milano\\Dipartimento di Matematica\\Via E. Bonardi, 9\\20133
Milano, Italy}
\email{irene.sabadini@polimi.it}

\title[A new resolvent equation for the S-functional calculus]
{A new resolvent equation for the S-functional calculus} \oddsidemargin
0.2in \evensidemargin 0.2in \topmargin -0.5in \textwidth
15.5truecm \textheight 23truecm

\def\R{\mathbb R}

\def\(s){\mathscr S(\R\times\R)}




\begin{document}
\maketitle
\tableofcontents
\parindent 0cm
\begin{abstract}
The S-functional calculus is a functional calculus for $(n+1)$-tuples of non necessarily commuting operators
 that can be considered a higher dimensional version of the classical Riesz-Dunford functional calculus for a single operator.  In this last calculus, the resolvent equation plays an important role in the proof of several results. Associated with the S-functional calculus there are two resolvent
  operators: the left  $S_L^{-1}(s,T)$ and the right one $S_R^{-1}(s,T)$, where $s=(s_0,s_1,\ldots ,s_n)\in \mathbb{R}^{n+1}$ and $T=(T_0,T_1,\ldots ,T_n)$ is an $(n+1)$-tuple of non commuting operators.
  These two S-resolvent operators satisfy the  S-resolvent equations
  $S_L^{-1}(s,T)s-TS_L^{-1}(s,T)=\mathcal{I}$,  and $sS_R^{-1}(s,T)-S_R^{-1}(s,T)T=\mathcal{I}$, respectively, where $\mathcal{I}$ denotes the identity operator. These equations allows to prove some properties of the S-functional calculus.
In this paper we  prove a new resolvent equation for the S-functional calculus which is the analogue of the classical resolvent equation. It is  interesting  to note
 that the equation  involves both the left and the right S-resolvent operators simultaneously.

\end{abstract}

\parindent 0cm

\section{Introduction}

The S-resolvent operators are a key tool in the definition of the higher dimensional version of the Riesz-Dunford functional calculus  called S-functional calculus.
This calculus works for $(n+1)$-tuples
$(T_0,T_1,\ldots ,T_n)$ of non necessarily commuting operators and  is based on the so-called $S$-spectrum, see \cite{DUKE, functionalcss}.
In the case of a single operator the $S$-functional calculus reduces to the Riesz-Dunford functional calculus (see \cite{ds, RS}).

When the operators $(T_0,T_1,\ldots ,T_n)$ commute among themselves, this calculus admits a commutative version called SC-functional calculus.
In this case the S-resolvent operator and the S-spectrum have a simpler expression, see \cite{SlaisComm}.

The class of functions on which this calculus is based is the so called set of slice hyperholomorphic (or slice monogenic) functions which are defined on subsets of
the Euclidean space $\mathbb{R}^{n+1}$ and have values in the Clifford algebra $\mathbb{R}_n$.

For more details on the S-functional calculus and the function theory on which it is based see the
 the monograph \cite{css_book}.

As it happens for the classical theory of monogenic functions (see \cite{BDS, csss, DSS, ghs}),
also in the class of slice hyperholomorphic functions
 there is the notion of left as well as of right hyperholomorphicity. But despite what happens in the monogenic case, for slice hyperholomorphic functions the Cauchy formulas for left and for right slice hyperholomorphic functions have two different kernels; moreover each of these kernels can be written in two different ways.

The calculus admits a quaternionic version,  which works for quaternionic linear operators and is based on slice hyperholomorphic (or slice regular) functions defined on subsets of the real algebra of quaternions $\mathbb{H}$ with values in the quaternions,  see \cite{cs, formulations}.
To explain our new result and its consequences, let us focus, at the moment, on the quaternionic setting which is simpler to illustrate.\\
Let us denote by $V$ a two sided quaternionic Banach space and
let $T:V\to V$ be a bounded right (or left) linear operator.
We recall that the  S-spectrum is defined as
$$
\sigma_S(T)=\{ s\in \mathbb{H}\ \ :\ \ T^2-2 Re(s)T+|s|^2\mathcal{I}\ \ \
{\rm is\ not\  invertible}\},
$$
where $s=s_0+s_1i+s_2j+s_3k$ is a quaternion, $Re(s)=s_0$,  $|s|^2=s_0^2+s_1^2+s_2^2+s_3^2$.
The left  and the right $S$-resolvent operators are defined as
\begin{equation}
S_L^{-1}(s,T):=-(T^2-2Re(s) T+|s|^2\mathcal{I})^{-1}(T-\overline{s}\mathcal{I}),\ \ \ s \in  \mathbb{H}\setminus\sigma_S(T)
\end{equation}
and
\begin{equation}
S_R^{-1}(s,T):=-(T-\overline{s}\mathcal{I})(T^2-2Re(s) T+|s|^2\mathcal{I})^{-1},\ \ \ s \in  \mathbb{H}\setminus\sigma_S(T),
\end{equation}
respectively.
The left $S$-resolvent
operator satisfies the equation
\begin{equation}\label{eSR1}
S_L^{-1}(s,T)s-TS_L^{-1}(s,T)=\mathcal{I},\ \ \ s \in  \mathbb{H}\setminus\sigma_S(T),
\end{equation}
and the right $S$-resolvent
operator satisfies
\begin{equation}\label{eSR2}
sS_R^{-1}(s,T)-S_R^{-1}(s,T)T=\mathcal{I},\ \ \ s \in  \mathbb{H}\setminus\sigma_S(T).
\end{equation}
Consider the complex plane $\mathbb{C}_I:=\mathbb{R}+I\mathbb{R}$, for $I\in \mathbb{S}$, where $\mathbb{S}$ is the unit sphere of purely imaginary quaternions. Observe that $\mathbb{C}_I$ can be identified with a complex plane since  $I^2=-1$ for every $I\in \mathbb{S}$.
Let $U\subset \hh$ be a suitable domain that contains the S-spectrum of $T$.
We define
  for  left slice hyperholomorphic functions $f:U \to \hh$ (see the precise definition in the sequel)  the quaternionic functional calculus as
\begin{equation}\label{quatinteg311def}
f(T)={{1}\over{2\pi }} \int_{\partial (U\cap \mathbb{C}_I)} S_L^{-1} (s,T)\  ds_I \ f(s),
\end{equation}
where $ds_I=- ds I$,
and for  right slice hyperholomorphic functions $f:U \to \hh$, we define
\begin{equation}\label{quatinteg311rightdef}
f(T)={{1}\over{2\pi }} \int_{\partial (U\cap \mathbb{C}_I)} \  f(s)\ ds_I \ S_R^{-1} (s,T).
\end{equation}
These definitions are well posed since the integrals do not depend neither on the open set $U$ and nor on the complex plane
$\mathbb{C}_I$ and can be extended to the case of $(n+1)$-tuples of operators, using slice hyperholomorphic functions with values in a Clifford algebra.
Using a similar notion of hyperholomorphicity and the S-spectrum in \cite{GMP} the authors introduce the
continuous functional calculus in a quaternionic Hilbert space.

The S-resolvent equations (\ref{eSR1}), (\ref{eSR2}) are useful to prove several properties of the S-functional calculus. However it is natural to ask if it is possible to obtain an analog of the classical resolvent equation
\begin{equation}\label{resolveqclassic}
(\lambda I-G)^{-1}(\mu I-G)^{-1}=\frac{(\lambda I-G)^{-1}-(\mu I-G)^{-1}}{\mu-\lambda},\ \ \lambda, \mu\in \mathbb{C}\setminus \sigma(G),
\end{equation}
where $G$ is a complex operator on a Banach space,
which might be useful to prove other properties of the calculus.
The main goal of this paper it to show that (\ref{resolveqclassic})
 can be generalized  in this non commutative setting,
 but it involves both the left and the right S-resolvent operators. Precisely, we will show that
\[
\begin{split}
S_R^{-1}(s,T)S_L^{-1}(p,T)&=[[S_R^{-1}(s,T)-S_L^{-1}(p,T)]p
\\
&
-\overline{s}[S_R^{-1}(s,T)-S_L^{-1}(p,T)]](p^2-2s_0p+|s|^2)^{-1},
\end{split}
\]
 for  $s$, $p \in  \mathbb{H}\setminus\sigma_S(T)$.
\\
It is also worthwhile to mention that
  the $S$-resolvent operator plays an important role in the
definition of the quaternionic version of the counterpart of the operator
$(I-zA)^{-1}$ in the realization $s(z)=D+zC(I-zA)^{-1}B$ for Schur multipliers, see \cite{acs1}.
The reader is referred to \cite{acs1, acs2, acs3} for Schur analysis in the slice hyperholomorphic setting
and to \cite{MR2002b:47144} and \cite{adrs} for an overview of Schur analysis in the complex setting.\\

It is interesting to note that in literature there are other cases in which the authors consider two resolvent operators.
We mention in particular the case of Schur analysis in the setting of upper
triangular operators and in the setting of compact Riemann surfaces. In the first case, the role of complex numbers is
played by diagonal operators and there are two ``point evaluations'' of an operator at a diagonal, one left and one right,
each corresponding to an associated resolvent operator; see \cite[(2.4)-(2.6), p. 256]{ALPZ}, but the resolvent equation
is related just with one resolvent at a time; see \cite[Corollary 2.9, p 266]{ALPZ}. In the setting of compact Riemann
surfaces (see \cite{vinnikov4,KLMV} for the general setting) there is a resolvent operator associated to every
meromorphic function on the given Riemann surface $X$
(see \cite[(4.1), p. 307]{av3}, and one needs two such operators, associated to a pair of functions which generate
the field of meromorphic functions on $X$, to study underlying spaces; see \cite[\S 5]{av3}. The same resolvent equation
is satisfied by all the resolvent operators; see \cite[Theorem 4.2, p. 309]{av3}.\\

In this setting both S-resolvent operators enter the resolvent equation.

The plan of the paper is as follows.

In Section 2 we recall some preliminary results on slice hyperholomorphic functions.

In Section 3
we state and prove the new resolvent equation and we show that there are two possible versions which are equivalent.
We prove our results for the S-functional calculus for $(n+1)$-tuples of non commuting operators and we show some applications of the resolvent equation.

In Section 4 we consider the  commutative version of the S-functional calculus, the so called SC-functional calculus
and we reformulate our main results for the quaternionic functional calculus. Since the proofs  follow the lines of the ones for the case of $(n+1)$-tuples of non commuting operators we will omit them in both cases.

\section{Preliminary results}
In this section we recall the notion of  slice hyperholomorphic functions and their Cauchy formulas,  see \cite{css_book}.

Let $\rr_n$ be the real Clifford algebra over $n$ imaginary units $e_1,\ldots ,e_n$
satisfying the relations $e_ie_j+e_je_i=0$,\  $i\not= j$, $e_i^2=-1.$
 An element in the Clifford algebra will be denoted by $\sum_A e_Ax_A$ where
$A=\{ i_1\ldots i_r\}\in \mathcal{P}\{1,2,\ldots, n\},\ \  i_1<\ldots <i_r$
 is a multi-index
and $e_A=e_{i_1} e_{i_2}\ldots e_{i_r}$, $e_\emptyset =1$.
An element $(x_0,x_1,\ldots,x_n)\in \mathbb{R}^{n+1}$  will be identified with the element
$
 x=x_0+\underline{x}=x_0+ \sum_{j=1}^nx_je_j\in\mathbb{R}_n
$
called paravector and the real part $x_0$ of $x$ will also be denoted by $Re(x)$.
The norm of $x\in\mathbb{R}^{n+1}$ is defined as $|x|^2=x_0^2+x_1^2+\ldots +x_n^2$.
 The conjugate of $x$ is defined by
$
\bar x=x_0-\underline x=x_0- \sum_{j=1}^nx_je_j.
$
Let
$$
\mathbb{S}=\{ \underline{x}=e_1x_1+\ldots +e_nx_n\ | \  x_1^2+\ldots +x_n^2=1\};
$$
for $I\in\mathbb{S}$ we obviously have $I^2=-1$.
Given an element $x=x_0+\underline{x}\in\rr^{n+1}$ let us set
$
I_x=\underline{x}/|\underline{x}|$ if $\underline{x}\not=0,
$
 and given an element $x\in\rr^{n+1}$, the set
$$
[x]:=\{y\in\rr^{n+1}\ :\ y=x_0+I |\underline{x}|, \ I\in \mathbb{S}\}
$$
is an $(n-1)$-dimensional sphere in $\mathbb{R}^{n+1}$.
The vector space $\mathbb{R}+I\mathbb{R}$ passing through $1$ and
$I\in \mathbb{S}$ will be denoted by $\mathbb{C}_I$ and
an element belonging to $\mathbb{C}_I$ will be indicated by $u+Iv$, for $u$, $v\in \mathbb{R}$.
With an abuse of notation we will write $x\in\mathbb{R}^{n+1}$.
Thus, if $U\subseteq\mathbb{R}^{n+1}$ is an open set,
a function $f:\ U\subseteq \mathbb{R}^{n+1}\to\mathbb{R}_n$ can be interpreted as
a function of the paravector $x$.

\begin{Dn}[Slice hyperholomorphic functions]
\label{defsmon}
Let $U\subseteq\mathbb{R}^{n+1}$ be an open set and let
$f: U\to\mathbb{R}_n$ be a real differentiable function. Let
$I\in\mathbb{S}$ and let $f_I$ be the restriction of $f$ to the
complex plane $\mathbb{C}_I$.
\\
The function  $f$ is said to be left slice  hyperholomorphic (or slice monogenic) if, for every
$I\in\mathbb{S}$, on $U\cap \mathbb{C}_I$ it satisfies
$$
\frac{1}{2}\left(\frac{\partial }{\partial u}f_I(u+Iv)+I\frac{\partial
}{\partial v}f_I(u+Iv)\right)=0.
$$
We will denote by $\mathcal{SM}(U)$ the set of left slice hyperholomorhic functions on the open set $U$ or by $\mathcal{SM}^L(U)$  when confusion may arise.
\\
The function $f$ is said to be right slice hyperholomorphic (or right slice monogenic) if,
for every
$I\in\mathbb{S}$, on $U\cap \mathbb{C}_I$, it satisfies
$$
\frac{1}{2}\left(\frac{\partial }{\partial u}f_I(u+Iv)+\frac{\partial
}{\partial v}f_I(u+Iv)I\right)=0.
$$
\\
We will denote by $\mathcal{SM}^R(U)$ the set of right slice hyperholomorphic functions on the open set $U$.
\end{Dn}

Slice hyperholomorphic functions possess good properties when they are defined on suitable domains which are
introduced in the following definition.

\begin{Dn}[Axially symmetric slice domain]\label{axsymm}
Let $U$ be a domain in $\rr^{n+1}$.
We say that $U$ is a
\textnormal{slice domain} (s-domain for short) if $U \cap \mathbb{R}$ is non empty and if $U\cap \mathbb{C}_I$ is a domain in $\mathbb{C}_I$ for all $I \in \mathbb{S}$.
We say that $U$ is
\textnormal{axially symmetric} if, for all $x \in U$, the
$(n-1)$-sphere $[x]$ is contained in $U$.
\end{Dn}

\begin{Dn}[Cauchy kernel for left slice hyperholomorphic functions]
Let $x$, $s\in \rr^{n+1}$ be such that $x\not\in [s]$.
Let $S_L^{-1}(s,x)$ be the function defined by
\begin{equation}\label{SL1}
S_L^{-1}(s,x):=-(x^2 -2x {\rm Re} [s]+|s|^2)^{-1}(x-\overline s).
\end{equation}
We say that  $S_L^{-1}(s,x)$ is the Cauchy kernel (for left slice hyperholomorphic functions) written in form I.
\end{Dn}
\begin{Pn}\label{ugualnuc}
Suppose that $x$ and $s\in\rr^{n+1}$ are such that $x\not\in [s]$.
The following identity holds:
\begin{equation}\label{third}
 -(x-\bar s)(x^2-2{\rm Re}(s)x+|s|^2)^{-1} = (s^2-2{\rm Re}(x)s+|x|^2)^{-1}(s-\bar x).
\end{equation}
\end{Pn}
\begin{Rk}{\rm By Proposition \ref{ugualnuc}
 $S_L^{-1}(s,x)$ can also be written as
\begin{equation}\label{SL2}
S_L^{-1}(s,x):=(s-\bar x)(s^2-2{\rm
Re}(x)s+|x|^2)^{-1}.
\end{equation}
In this case, we will say  $S_L^{-1}(s,x)$ is written in form II.
}
\end{Rk}
\begin{Pn}\label{pro14}
The function $S_L^{-1}(s,x)$ is left slice hyperholomorphic in the variable $x$
and right slice hyperholomorphic in the variable $s$  for $x\not\in [s]$.
\end{Pn}
The case of the Cauchy kernel for right slice hyperholomorphic functions is similar.
\begin{Dn}[Cauchy kernel for right slice hyperholomorphic functions]\label{CauchykernelrightMON}
Let $x$, $s\in \rr^{n+1}$ be such that $x\not\in [s]$.
The Cauchy kernel  $S_R^{-1}(s,x)$ for right slice hyperholomorphic functions is defined by
\begin{equation}\label{SR1}
S_R^{-1}(s,x):= -(x-\bar s)(x^2-2{\rm Re}(s)x+|s|^2)^{-1}.
\end{equation}
We say that  $S_R^{-1}(s,x)$ is written in form I.
\end{Dn}
\begin{Rk}{\rm
An analog of Proposition \ref{ugualnuc} holds in fact:
\begin{equation}\label{third}
-(x-\bar s)(x^2-2{\rm Re}(s)x+|s|^2)^{-1}=(s^2-2{\rm Re}(x)s+|x|^2)^{-1}(s-\bar x),
\end{equation}
for $x,$ $s\in\rr^{n+1}$ such that $x\not\in [s]$.\\
Thus $S_R^{-1}(s,x)$ can be written as
$$
S_R^{-1}(s,x)= (s^2-2{\rm Re}(x)s+|x|^2)^{-1}(s-\bar x),
$$
and in this case we say that $S_R^{-1}(s,x)$ is written in form II.
}
\end{Rk}
\begin{Tm}[The Cauchy formula with slice hyperholomorphic kernel]
\label{Cauchygenerale}
Let $U\subset\mathbb{R}^{n+1}$ be an axially symmetric s-domain.
Suppose that $\partial (U\cap \mathbb{C}_I)$ is a finite union of
continuously differentiable Jordan curves  for every $I\in\mathbb{S}$.  Set  $ds_I=-ds I$ for $I\in \mathbb{S}$.
\begin{itemize}
\item
If $f$ is
a (left) slice hyperholomorphic function on a set that contains $\overline{U}$ then
\begin{equation}\label{Cauchyleft}
 f(x)=\frac{1}{2 \pi}\int_{\partial (U\cap \mathbb{C}_I)} S_L^{-1}(s,x)ds_I f(s)
\end{equation}
 and the
integral does not depend on $U$ and on the imaginary unit
$I\in\mathbb{S}$.
\item
If $f$ is a right slice hyperholomorphic function on a set that contains $\overline{U}$,
then
\begin{equation}\label{Cauchyright}
 f(x)=\frac{1}{2 \pi}\int_{\partial (U\cap \mathbb{C}_I)}  f(s)ds_I S_R^{-1}(s,x)
 \end{equation}
and the integral   does not depend on  $U$ and on  the imaginary unit $I\in\mathbb{S}$.
\end{itemize}
\end{Tm}
The above Cauchy formulas are the  starting point to define the S-functional calculus.
A crucial fact of slice hyperholomorphic functions in the representation formula (also called  structure formula).
This formula will be used in the sequel to give applications of the new resolvent equation.

\begin{Tm}[Representation Formula]\label{formulaMON} Let $U$ be an axially symmetric s-domain $U \subseteq  \mathbb{H}$.
\begin{itemize}
\item
Let $f$ be a (left) slice hyperholomorphic function on $U$.  Choose any
$J\in \mathbb{S}$.  Then the following equality holds for all $x=u+Iv \in U $:
\begin{equation}\label{distribution}
f(u+Iv) =\frac{1}{2}\Big[   f(u+Jv)+f(u-Jv)\Big] +I\frac{1}{2}\Big[ J[f(u-Jv)-f(u+Jv)]\Big].
\end{equation}
Moreover, for all $u, v \in \mathbb{R}$ such that $u+v\mathbb{S} \subseteq U $, there exist $\mathbb{R}_n$-valued functions $\alpha$, $\beta$ depending on $u,v$ only such that for all $K \in \mathbb{S}$
\begin{equation}\label{cappa}
\frac{1}{2}\Big[   f(u+Kv)+f(u-Kv)\Big]=\alpha(u,v) \ \ {\sl and}
\ \  \frac{1}{2}\Big[ K[f(u-Kv)-f(u+Kv)]\Big]=\beta(u,v).
\end{equation}
\item
Let $f$ be a right slice hyperholomorphic function on $U$.  Choose any
$J\in \mathbb{S}$.  Then the following equality holds for all $x=u+Iv \in U $:
\begin{equation}\label{distributionright}
f(u+Iv) =\frac{1}{2}\Big[   f(u+Jv)+f(u-Jv)\Big]+\frac{1}{2}\Big[ [f(u-Jv)-f(u+Jv)]J\Big]I.
\end{equation}
Moreover, for all $u, v \in \mathbb{R}$ such that $u+v\mathbb{S} \subseteq U $, there exist $\mathbb{R}_n$-valued functions $\alpha$, $\beta$ depending on $u,v$ only such that for all $K \in \mathbb{S}$
\begin{equation}\label{capparight}
\frac{1}{2}\Big[   f(u+Kv)+f(u-Kv)\Big]=\alpha(u,v) \ \ {\sl and}
\ \ \frac{1}{2}\Big[ [f(u-Kv)-f(u+Kv)]K\Big]=\beta(u,v).
\end{equation}
\end{itemize}
\end{Tm}

\section{The case of several non commuting operators}

In the sequel, we will consider a Banach space $V$ over
$\mathbb{R}$
 with norm $\|\cdot \|$.
It is possible to endow $V$
with an operation of multiplication by elements of $\rr_n$ which gives
a two-sided module over $\rr_n$.
A two-sided module $V$ over $\rr_n$ is called a Banach module over $\rr_n$,
 if there exists a constant $C \geq 1$  such
that $\|va\|\leq C\| v\| |a|$ and $\|av\|\leq C |a|\| v\|$ for all
$v\in V$ and $a\in\rr_n$.
 By $V_n$ we denote $V\otimes \rr_n$ over $\rr_n$; $V_n$ turns out to be a
 two-sided Banach module .\\
 An element in $V_n$ is of the type $\sum_A v_A\otimes e_A$ (where
 $A=i_1\ldots i_r$, $i_\ell\in \{1,2,\ldots, n\}$, $i_1<\ldots <i_r$ is a multi-index).
The multiplications of an element $v\in V_n$ with a scalar
$a\in \rr_n$ are defined by $va=\sum_A v_A \otimes (e_A a)$ and $av=\sum_A v_A \otimes (ae_A )$.
For simplicity, we will write
$\sum_A v_A e_A$ instead of $\sum_A v_A \otimes e_A$. Finally, we define $\| v\|^2_{V_n}=
\sum_A\| v_A\|^2_V$.

We denote by
$\mathcal{B}(V)$  the space
of bounded $\mathbb{R}$-homomorphisms of the Banach space $V$ to itself
 endowed with the natural norm denoted by $\|\cdot\|_{\mathcal{B}(V)}$.
Given $T_A\in \mathcal{B}(V)$, we can introduce the operator $T=\sum_A T_Ae_A$ and
its action on $v=\sum v_Be_B\in V_n$ as $T(v)=\sum_{A,B}
T_A(v_B)e_Ae_B$. The operator $T$ is a right-module homomorphism which is a bounded linear
map on $V_n$.
\\
In the sequel, we will consider operators of the form
$T=T_0+\sum_{j=1}^ne_jT_j$ where $T_j\in\mathcal{B}(V)$ for $j=0,1,\ldots ,n$.
The subset of such operators in ${\mathcal{B}(V_n)}$ will be denoted by $\mathcal{B}^{\small 0,1}(V_n)$.
We define $\|T\|_{\mathcal{B}^{\small 0,1}(V_n)}=\sum_j \|T_j\|_{\mathcal{B}(V)}$.
Note that, in the sequel, we will omit the subscript $\mathcal{B}^{\small 0,1}(V_n)$ in the norm of an operator. Note also that  $\|TS\|\leq \|T\| \|S\|$.

\begin{Dn}
Let $T\in\mathcal{B}^{\small 0,1}(V_n)$. We define the left Cauchy kernel operator series
 or $S$-resolvent operator series  as
\begin{equation}\label{quatCKSLoperator}
S_L^{-1}(s,T)=\sum_{n\geq 0} T^n s^{-1-n},
\end{equation}
and the right Cauchy kernel operator series as
\begin{equation}\label{quatCKSRightoperator}
S_R^{-1}(s,T)=\sum_{n\geq 0} s^{-1-n}T^n ,
\end{equation}
for $\|T\| <|s|$.
\end{Dn}
The Cauchy kernel operator series are the power series expansion of the S-resolvent operators. Their sum is computed in the following result:
\begin{Tm}\label{Ssinistro}
Let $T\in\mathcal{B}^{\small 0,1}(V_n)$ and let $s \in \mathbb{H}$.
Then,  for $\|T\|< |s|$, we have
\begin{equation}\label{SERLEFTgf}
\sum_{m\geq 0} T^m s^{-1-m}=-(T^2-2Re(s) T+|s|^2\mathcal{I})^{-1}(T-\overline{s}\mathcal{I}),
\end{equation}
\begin{equation}\label{SERRIGHTgf}
\sum_{m\geq 0} s^{-1-m} T^m =-(T-\overline{s}\mathcal{I})(T^2-2Re(s) T+|s|^2\mathcal{I})^{-1}.
\end{equation}
\end{Tm}
We observe that the sum of the above series are independent of the fact that the components of the paravector operator $T$ commute. Moreover the operators on right hand sides of
(\ref{SERLEFTgf}) and (\ref{SERRIGHTgf}) are defined on a subset of $\mathbb{R}^{n+1}$
that is larger then $\{ s\in \mathbb{R}^{n+1} \ :\ \|T\|< |s|\}.$ This fact suggests the definition of $S$-spectrum, of $S$-resolvent set and of $S$-resolvent operators.

\begin{Dn}[The $S$-spectrum and the $S$-resolvent set]
Let $T\in\mathcal{B}^{\small 0,1}(V_n)$.
We define the $S$-spectrum $\sigma_S(T)$ of $T$  as:
$$
\sigma_S(T)=\{ s\in \mathbb{R}^{n+1}\ \ :\ \ T^2-2 \ Re (s)T+|s|^2\mathcal{I}\ \ \
{\it is\ not\  invertible}\}.
$$
The $S$-resolvent set $\rho_S(T)$ is defined by
$$
\rho_S(T)=\mathbb{R}^{n+1}\setminus\sigma_S(T).
$$
\end{Dn}

\begin{Dn}[The $S$-resolvent operators]
Let  $T\in \mathcal{B}^{\small 0,1}(V_n)$ and $s\in \rho_S(T)$.
We define the left $S$-resolvent operator as
\begin{equation}\label{quatSresolrddlft}
S_L^{-1}(s,T):=-(T^2-2Re(s) T+|s|^2\mathcal{I})^{-1}(T-\overline{s}\mathcal{I}),
\end{equation}
and the right $S$-resolvent operator as
\begin{equation}\label{quatSresorig}
S_R^{-1}(s,T):=-(T-\overline{s}\mathcal{I})(T^2-2Re(s) T+|s|^2\mathcal{I})^{-1}.
\end{equation}
\end{Dn}
The operators $S_L^{-1}(s,T)$ and $S_R^{-1}(s,T)$ satisfy the equations below, see \cite{css_book}:
\begin{Tm}
Let $T\in \mathcal{B}^{\small 0,1}(V_n)$ and let $s \in \rho_S(T)$. Then, the left $S$-resolvent
operator satisfies the equation
\begin{equation}\label{quatSresolrddlftequ}
S_L^{-1}(s,T)s-TS_L^{-1}(s,T)=\mathcal{I},
\end{equation}
and the right $S$-resolvent
operator satisfies the equation
\begin{equation}\label{quatSresorigequa}
sS_R^{-1}(s,T)-S_R^{-1}(s,T)T=\mathcal{I}.
\end{equation}
\end{Tm}

 Our goal is to establish the analogue of the classical resolvent equation. To this end, we need some preliminary results. A crucial fact is the following Theorem \ref{Ssinistro} that will give us the hint to discover what is the structure of the resolvent equation in this non commutative setting at least in the case the S-resolvent operators are expressed in power series.
\begin{Tm}\label{Ssinistro}
Let $A$, $B\in\mathcal{B}(V_n)$ and let $s$, $p \in \mathbb{R}^{n+1}$.
Then,  for $|p|< |s|$, we have
\begin{equation}\label{SERLEFT}
\sum_{m\geq 0} p^m A s^{-1-m}=-(p^2-2Re(s) p+|s|^2)^{-1}(pA-A\overline{s}),
\end{equation}
and
\begin{equation}\label{SERRIGHT}
\sum_{m\geq 0} s^{-1-m} Bp^m =-(Bp-\overline{s}B)(p^2-2Re(s) p+|s|^2)^{-1}.
\end{equation}
Moreover, (\ref{SERRIGHT}) can be written as
\begin{equation}\label{SERRIGHTbis}
\sum_{m\geq 0} s^{-1-m} Bp^m =(s^2-2Re(p) s+|p|^2)^{-1}(sB-B\overline{p}).
\end{equation}
\end{Tm}
\begin{proof}
To verify (\ref{SERLEFT}) define
$$
X:=(p^2-2Re(s) p+|s|^2)\sum_{m\geq 0} p^m A s^{-1-m}
$$
and observe that
\begin{equation}
\begin{split}
X&=\sum_{m\geq 0} (p^2-2Re(s) p+|s|^2)p^m A s^{-1-m}
\\
&
= p^2A s^{-1}-2Re(s) pA s^{-1}+|s|^2 A s^{-1}
\\
&
+p^3 A s^{-2}-2Re(s) p^2 A s^{-2}+|s|^2p A s^{-2}
\\
&
+p^4 A s^{-3}-2Re(s) p^3 A s^{-3}+|s|^2p^2 A s^{-3}+\ldots
\\
&
=-(pA-A\overline{s})+\sum_{m\geq 2} p^m A (s^2-2Re(s)s+|s|^2) s^{-1-m}.
\end{split}
\end{equation}
Since any paravector $s$ satisfies
$$
s^2-2Re(s)s+|s|^2=0
$$
we deduce that
$$
X=(p^2-2Re(s) p+|s|^2)\sum_{m\geq 0} p^m A s^{-1-m}=-(pA-A\overline{s})
$$
and the statement follows.
The  equality in (\ref{SERRIGHT}) can be verified by setting
$$
Y:=\sum_{m\geq 0} s^{-1-m} Bp^m (p^2-2Re(s) p+|s|^2)
$$
and observing that
$$
Y=-(Bp-\overline{s}B)+\sum_{m\geq 0} s^{-1-m} Bp^m (p^2-2Re(p) p+|p|^2)= -(Bp-\overline{s}B).
$$
With similar computations one can verify equality (\ref{SERRIGHTbis}).
\end{proof}
\begin{Cy}\label{asdf}
Let $A$, $B\in\mathcal{B}(V_n)$ and let $s$, $p$ be paravectors.
Then,  for $|p|< |s|$, the following equations hold
\begin{equation}\label{SERLEFTABfddf}
\begin{split}
\sum_{j=0}^m p^j A s^{-1-j}&=-(p^2-2Re(s) p+|s|^2)^{-1}(pA-A\overline{s})
\\
&
+p^{m+1}(p^2-2Re(s) p+|s|^2)^{-1}(pA-A\overline{s})s^{-1-m},
\end{split}
\end{equation}
and
\begin{equation}\label{SERRIGHTAB}
\begin{split}
\sum_{j=0}^m s^{-1-j} Bp^j &=-(Bp-\overline{s}B)(p^2-2Re(s) p+|s|^2)^{-1}
\\
&+
s^{-1-m}(Bp-\overline{s}B)(p^2-2Re(s) p+|s|^2)^{-1}p^{m+1}.
\end{split}
\end{equation}
Moreover, (\ref{SERRIGHTAB}) can also written as
\begin{equation}\label{SERRIGHTABbis}
\begin{split}
\sum_{j=0}^m s^{-1-j} Bp^j &=(s^2-2Re(p) s+|p|^2)^{-1}(sB-B\overline{p})
\\
&-
s^{-1-m}(s^2-2Re(p) s+|p|^2)^{-1}(sB-B\overline{p})p^{m+1}.
\end{split}
\end{equation}
\end{Cy}
\begin{proof}
Identity (\ref{SERLEFTABfddf}) follows from
\[
\sum_{j=0}^m p^j A s^{-1-j}=\sum_{j=0}^\infty p^j A s^{-1-j}
-\sum_{j=m+1}^\infty p^j A s^{-1-j},
\]
that can be written as
\[
\sum_{j=0}^m p^j A s^{-1-j}=\sum_{j=0}^\infty p^j A s^{-1-j}
-p^{m+1}\Big(\sum_{j=0}^\infty p^j A s^{-1-j}\Big)s^{-1-m},
\]
but now we use (\ref{SERLEFT}) to get the result.
Identity (\ref{SERRIGHTAB}) and (\ref{SERRIGHTABbis}) follow with similar computations.
\end{proof}
We now prove the new S-resolvent equation. In the proof we first consider
the case in which the S-resolvent operators admit the power series expansion
$$
S_L^{-1}(s,T)=\sum_{m\geq 0} T^m s^{-1-m},
\ \ \ \ \ \
S_R^{-1}(s,T)=\sum_{m\geq 0} s^{-1-m}T^m,
$$
that is for $\|T\| <|s|$.
Then we verify that such equation holds  in general.

\begin{Tm}\label{RLRESOLVEQ}
Let $T\in\mathcal{B}^{\small 0,1}(V_n)$ and let $s$ and $p\in \rho_S(T)$. Then we have
\begin{equation}\label{RLresolv}\small
S_R^{-1}(s,T)S_L^{-1}(p,T)=((S_R^{-1}(s,T)-S_L^{-1}(p,T))p
-\overline{s}(S_R^{-1}(s,T)-S_L^{-1}(p,T)))(p^2-2s_0p+|s|^2)^{-1}.
\end{equation}
Moreover, the resolvent equation can also be written as
\begin{equation}\label{RLresolvII}\small
S_R^{-1}(s,T)S_L^{-1}(p,T)=(s^2-2p_0s+|p|^2)^{-1}(s(S_R^{-1}(s,T)-S_L^{-1}(p,T))
-(S_R^{-1}(s,T)-S_L^{-1}(p,T))\overline{p} ).
\end{equation}
\end{Tm}
\begin{proof} We prove the theorem in two steps.
\\
STEP I.
First we assume that the $S$-resolvent operators are expressed in power series. If $\|T\|<|p| <|s|$  then
the S-resolvent operators have power series expansion and so
\begin{equation}\label{seriesresolvent}
S_R^{-1}(s,T) S_L^{-1}(p,T)=(\sum_{j\geq 0} s^{-1-j} T^j)(\sum_{j\geq 0} T^jp^{-1-j} ).
\end{equation}
By setting
$$
\Lambda_m(s,p;T):=\sum_{j=0}^ms^{-1-j}(T^mp^{-1-m})p^{j}
$$
(\ref{seriesresolvent}) can be written as
$$
S_R^{-1}(s,T) S_L^{-1}(p,T)=\sum_{m\geq 0}\Lambda_m(s,p;T).
$$
Formula (\ref{SERRIGHTAB}) with $B=T^mp^{-1-m}$ and some computations give
\begin{equation}
\begin{split}
\Lambda_m(s,p;T)&=
-((T^mp^{-1-m})p-\overline{s}(T^mp^{-1-m}))(p^2-2Re(s) p+|s|^2)^{-1}
\\
&+
s^{-1-m}((T^mp^{-1-m})p-\overline{s}(T^mp^{-1-m}))(p^2-2Re(s) p+|s|^2)^{-1}p^{m+1}\\
&=
-[(T^mp^{-1-m})p-\overline{s}(T^mp^{-1-m})
\\
&+
(s^{-1-m}T^m)p-\overline{s}(s^{-1-m}T^m)](p^2-2Re(s) p+|s|^2)^{-1}.
\end{split}
\end{equation}
From the chain of equalities
\begin{equation}
\begin{split}
S_R^{-1}(s,T) S_L^{-1}(p,T)&=\sum_{m\geq 0}\Lambda_m(s,p;T)
\\
&=
-[(\sum_{m\geq 0}(T^mp^{-1-m})p-\overline{s}\sum_{m\geq 0}(T^mp^{-1-m}))]
\\
&+
(\sum_{m\geq 0}s^{-1-m}T^m)p-\overline{s}\sum_{m\geq 0}s^{-1-m}T^m)](p^2-2Re(s) p+|s|^2)^{-1}
\end{split}
\end{equation}
(\ref{RLresolv}) follows.

To prove that the resolvent equation can be written in the second form (\ref{RLresolvII}) observe that $\Lambda_m(s,p;T)$ can also be written using (\ref{SERRIGHTbis}) as
\begin{equation}
\begin{split}
\Lambda_m(s,p;T)&=(s^2-2Re(p) s+|p|^2)^{-1}(s(T^mp^{-1-m})-(T^mp^{-1-m})\overline{p})
\\
&-
s^{-1-m}(s^2-2Re(p) s+|p|^2)^{-1}(s(T^mp^{-1-m})-(T^mp^{-1-m})\overline{p})p^{m+1}.
\end{split}
\end{equation}
so taking the sum $\sum_{m\geq 0}\Lambda_m(s,p;T)$ we get the second version of the resolvent equation.
\\
\\
STEP II.
We prove that,
for $s$ and $p\in \rho_S(T)$,
(\ref{RLresolv}) and (\ref{RLresolvII}) hold with $S_R^{-1}(s,T)$ and  $S_L^{-1}(p,T)$ defined in (\ref{quatSresolrddlft}) and (\ref{quatSresorig}), respectively.
\\
Let us verify  (\ref{RLresolv}).
Since  $s$ and $p\in \rho_S(T)$ the left and right S-resolvent operators defined by
(\ref{quatSresolrddlft}) and (\ref{quatSresorig}) satisfy the left and the right resolvent equations
 (\ref{quatSresolrddlftequ}) and (\ref{quatSresorigequa}), respectively.
To verify (\ref{RLresolv}) we have to show that
$S_R^{-1}(s,T)S_L^{-1}(p,T)(p^2-2s_0p+|s|^2)$ equals
$$
(S_R^{-1}(s,T)-S_L^{-1}(p,T))p-\overline{s}(S_R^{-1}(s,T)-S_L^{-1}(p,T)).
$$
 To do this we use the left and the right S-resolvent equations \eqref{quatSresolrddlftequ},
\eqref{quatSresorigequa}.
Indeed, using the left S-resolvent equation, written as
$$
S_L^{-1}(p,T)p=TS_L^{-1}(p,T)+\mathcal{I},
$$
 we have
\begin{equation}
\begin{split}
& S_R^{-1}(s,T)S_L^{-1}(p,T)(p^2-2s_0p+|s|^2)
=S_R^{-1}(s,T)[S_L^{-1}(p,T)p]p
\\
&
-2s_0S_R^{-1}(s,T)S_L^{-1}(p,T)p
+|s|^2S_R^{-1}(s,T)S_L^{-1}(p,T)
\\
&
=S_R^{-1}(s,T) [TS_L^{-1}(p,T)+\mathcal{I}]p
-2s_0S_R^{-1}(s,T)[TS_L^{-1}(p,T)+\mathcal{I}]
\\
&
+|s|^2S_R^{-1}(s,T)S_L^{-1}(p,T)
\end{split}
\end{equation}
and using again the left $S$-resolvent equation
\begin{equation}
\begin{split}
&S_R^{-1}(s,T)S_L^{-1}(p,T)(p^2-2s_0p+|s|^2)
=S_R^{-1}(s,T) T[TS_L^{-1}(p,T)+\mathcal{I}]+S_R^{-1}(s,T)p
\\
&
-2s_0S_R^{-1}(s,T)[TS_L^{-1}(p,T)+\mathcal{I}]
+|s|^2S_R^{-1}(s,T)S_L^{-1}(p,T)
\end{split}
\end{equation}
we obtain
\begin{equation}
\begin{split}
S_R^{-1}(s,T)&S_L^{-1}(p,T)(p^2-2s_0p+|s|^2)
\\
&
=[S_R^{-1}(s,T)T] TS_L^{-1}(p,T)+S_R^{-1}(s,T) T]+S_R^{-1}(s,T)p
\\
&
-2s_0[[S_R^{-1}(s,T)T]S_L^{-1}(p,T)+S_R^{-1}(s,T)]
\\
&
+|s|^2S_R^{-1}(s,T)S_L^{-1}(p,T).
\end{split}
\end{equation}
Now we use the right $S$-resolvent equation
$$
S_R^{-1}(s,T)T=sS_R^{-1}(s,T)-\mathcal{I}
$$
we obtain
\begin{equation}
\begin{split}
S_R^{-1}(s,T)&S_L^{-1}(p,T)(p^2-2s_0p+|s|^2)
\\
&
=[[sS_R^{-1}(s,T)-\mathcal{I}] T]S_L^{-1}(p,T)+sS_R^{-1}(s,T)-\mathcal{I}]+S_R^{-1}(s,T)p
\\
&
-2s_0[[sS_R^{-1}(s,T)-\mathcal{I}]S_L^{-1}(p,T)+S_R^{-1}(s,T)]
\\
&
+|s|^2S_R^{-1}(s,T)S_L^{-1}(p,T).
\end{split}
\end{equation}
Iterating the use of the above right $S$-resolvent equation we get
\begin{equation}
\begin{split}
S_R^{-1}(s,T)&S_L^{-1}(p,T)(p^2-2s_0p+|s|^2)
\\
&
=[ s[sS_R^{-1}(s,T)-\mathcal{I}]-T]S_L^{-1}(p,T)+sS_R^{-1}(s,T)-\mathcal{I}]+S_R^{-1}(s,T)p
\\
&
-2s_0[[sS_R^{-1}(s,T)S_L^{-1}(p,T)-S_L^{-1}(p,T)]+S_R^{-1}(s,T)]
\\
&
+|s|^2S_R^{-1}(s,T)S_L^{-1}(p,T),
\end{split}
\end{equation}
which leads to
\begin{equation}
\begin{split}
S_R^{-1}(s,T)&S_L^{-1}(p,T)(p^2-2s_0p+|s|^2)
\\
&
=(s^2-2s_0s+|s|^2)S_R^{-1}(s,T)S_L^{-1}(p,T)
\\
&
+[S_R^{-1}(s,T)-S_L^{-1}(p,T)]p-\overline{s}[S_R^{-1}(s,T)-S_L^{-1}(p,T)],
\end{split}
\end{equation}
and since $s^2-2s_0s+|s|^2=0$ we obtain (\ref{RLresolv}).
With similar computations we can show that also (\ref{RLresolvII}) holds.
\end{proof}
We now observe that in the commutative case besides the resolvent equation, also the following relation between the resolvent operators
$$
(\lambda I-G)^{-1}(\mu I-G)^{-1}=(\mu I-G)^{-1}(\lambda I-G)^{-1}, \ \ \ {\rm for} \ \ \ \lambda, \mu \in \rho(G)
$$
holds.
In the non commutative case  we cannot aspect the validity of such a relation, however we will show that
an analogous equation holds for the so-called pseudo S-resolvent operators defined below.
\begin{Dn} Let $T\in \mathcal{B}^{0,1}(V_n)$. We define, for $s\in \rho_S(T)$, the pseudo S-resolvent operator of $T$ is defined as
$$
Q_s(T):=(T^2-2Re(s) T+|s|^2\mathcal{I})^{-1}.
$$
\end{Dn}
With the above definition the resolvents $S_L^{-1}(s,T)$ and $S_R^{-1}(s,T)$ become
\begin{equation}
S_L^{-1}(s,T):=-Q_s(T)(T-\overline{s}\mathcal{I}),\ \ \ s \in  \rho_S(T),
\end{equation}
and
\begin{equation}
S_R^{-1}(s,T):=-(T-\overline{s}\mathcal{I})Q_s(T),\ \ \ s \in  \rho_S(T).
\end{equation}
We now prove the following:
\begin{Tm}
Let $T\in \mathcal{B}^{0,1}(V_n)$ and let $s$, $p\in \rho_S(T)$. Then we have
$$
(T-\overline{s}\mathcal{I})Q_s(T)Q_p(T)(T-\overline{p}\mathcal{I})
=(T-\overline{s}\mathcal{I})Q_p(T)Q_s(T)(T-\overline{p}\mathcal{I}).
$$
\end{Tm}
\begin{proof}
It follows from the fact that
\begin{equation}
\begin{split}
(T^2-2Re(s) T+|s|^2\mathcal{I})&(T^2-2Re(p) T+|p|^2\mathcal{I})
\\
&=
(T^2-2Re(p) T+|p|^2\mathcal{I})(T^2-2Re(s) T+|s|^2\mathcal{I}).
\end{split}
\end{equation}
Since  $s$, $p\in \rho_S(T)$ we can take the inverse and the statement follows.
\end{proof}

\begin{Rk}{\rm
Observe that the function $F_T(s,p)$ defined by
$$
F_T(s,p):=S_R^{-1}(s,T)S_L^{-1}(p,T)
$$
is left slice hyperholomorphic in $s$ and it is right slice hyperholomorphic in $p$ with values in $\mathcal{B}(V_n)$.
The function
$$
G_T(s,p):=S_L^{-1}(p,T)S_R^{-1}(s,T)
$$
is not slice hyperholomorphic neither in $p$ nor in $s$.
}
\end{Rk}

\begin{Rk}{\rm Using the star products left and right in the variables $s,p$, which will be denoted by $\star_{s, left}$, $\star_{p, right}$ respectively, see \cite{acls}, the resolvent equation (\ref{RLresolv}) can be written as
$$
S_R^{-1}(s,T)S_L^{-1}(p,T)= [S_R^{-1}(s,T)-S_L^{-1}(p,T)]\star_{s, left} (p-\overline{s})(p^2-2Re(s) p+|s|^2)^{-1}\mathcal{I},
$$
or
$$
S_R^{-1}(s,T)S_L^{-1}(p,T)= (s-\overline{p})(s^2-2Re(p) s+|p|^2)^{-1}\mathcal{I}\star_{p, right}[S_R^{-1}(s,T)-S_L^{-1}(p,T)].
$$
}
\end{Rk}

\subsection{Some applications}
\medskip
Here we recall the formulations of the S-functional calculus and the we use
 the resolvent equation to deduce some results.

We first recall two important properties of
the S-spectrum.
\begin{Tm}[Structure of the $S$-spectrum]\label{strutturaS}
\par\noindent
Let $T\in\mathcal{B}^{\small 0,1}(V_n)$
and suppose that $p=p_0+\underline{p}$ belongs  $\sigma_S(T)$
 with $ \underline{p}\neq 0$.
Then all the elements of the $(n-1)$-sphere $[p]$ belong to $\sigma_S(T)$.
\end{Tm}
This result implies that if $p\in\sigma_S(T)$ then either $p$ is a
real point or the whole $(n-1)$-sphere $[p]$ belongs to
$\sigma_S(T)$.

\par\noindent
\begin{Tm}[Compactness of $S$-spectrum]\label{compattezaS}
 Let $T\in\mathcal{B}^{\small 0,1}(V_n)$. Then
the $S$-spectrum $\sigma_S (T)$  is a compact nonempty set.
Moreover, $\sigma_S (T)$ is
contained in $\{s\in\rr^{n+1}\, :\,  |s|\leq \|T\| \ \}$.
\end{Tm}

\begin{Dn}
  Let $V_n$ be a two sided Banach module,   $T\in\mathcal{B}^{0,1}(V_n)$and let $U \subset \mathbb{R}^{n+1}$ be an axially symmetric s-domain
that contains  the $S$-spectrum $\sigma_S(T)$  such that
$\partial (U\cap \mathbb{C}_I)$ is the union of a finite number of
continuously differentiable Jordan curves  for every $I\in\mathbb{S}$.
In this case we say that $U$ is a $T$-admissible open set.
\end{Dn}
We can now introduce the class of functions for which we can define the two versions of the S-functional calculus.
 \begin{Dn}\label{quatdef3.9}
Let $V_n$ be a two sided Banach module,   $T\in\mathcal{B}^{0,1}(V_n)$  and let  $W$ be an open set in $\mathbb{R}^{n+1}$.
\begin{itemize}
\item[(i)]
A function  $f\in \mathcal{SM}^L(W)$  is said to be locally left hyperholomorphic  on $\sigma_S(T)$
if there exists a $T$-admissible domain $U\subset \mathbb{R}^{n+1}$ such that $\overline{U}\subset W$, on
which $f$ is left slice hyperholomorphic.
We will denote by $\mathcal{SM}^L_{\sigma_S(T)}$ the set of locally
left hyperholomorphic functions on $\sigma_S (T)$.
\item[(ii)]
A function $f\in \mathcal{SM}^R(W)$ is said to be locally right regular on $\sigma_S(T)$
if there exists a $T$-admissible domain $U\subset \mathbb{R}^{n+1}$ such that $\overline{U}\subset W$, on
which $f$ is right slice hyperholomorphic.
We will denote by $\mathcal{SM}^R_{\sigma_S(T)}$ the set of locally
right slice hyperholomorphic functions on $\sigma_S (T)$.
\end{itemize}
\end{Dn}

\begin{Dn}[The  $S$-functional calculus]
Let $V_n$ be a two sided  Banach module and  $T\in\mathcal{B}^{0,1}(V_n)$.
  Let $U\subset  \mathbb{R}^{n+1}$ be a $T$-admissible domain and set $ds_I=- ds I$. We define
\begin{equation}\label{Scalleft}
f(T)={{1}\over{2\pi }} \int_{\partial (U\cap \mathbb{C}_I)} S_L^{-1} (s,T)\  ds_I \ f(s), \ \ {\it for}\ \ f\in \mathcal{SM}^L_{\sigma_S(T)},
\end{equation}
and
\begin{equation}\label{Scalright}
f(T)={{1}\over{2\pi }} \int_{\partial (U\cap \mathbb{C}_I)} \  f(s)\ ds_I \ S_R^{-1} (s,T),\ \  {\it for}\ \ f\in \mathcal{SM}^R_{\sigma_S(T)}.
\end{equation}
\end{Dn}

We now define the Riesz projectors for the S-functional calculus.
We begin with a preliminary lemma.

\begin{La}\label{Lemma321gf}
Let $B\in \mathcal{B}(V_n)$ and let $G$ be an axially symmetric s-domain such that  $p\in G$.
Then
\begin{equation}\label{efdhsslj}
(\overline{s}B-Bp)(p^2-2s_0p+|s|^2)^{-1}=
(s^2-2p_0s+|p|^2)^{-1}(sB-B\overline{p}), \ \ \ \  p\not\in [s],
\end{equation}
and
\begin{equation}\label{efdh}
\frac{1}{2\pi}\int_{\partial(G\cap\mathbb{C}_I)}ds_I
(\overline{s}B-Bp)(p^2-2s_0p+|s|^2)^{-1}=B.
\end{equation}
\end{La}
\begin{proof}
Formula (\ref{efdhsslj}) is obtained by direct computation.
Let us prove (\ref{efdh}),
so we write
\[
\begin{split}
\frac{1}{2\pi}\int_{\partial(G\cap\mathbb{C}_I)}&ds_I
(\overline{s}B-Bp)(p^2-2s_0p+|s|^2)^{-1}
\\
&
=
\frac{1}{2\pi}\int_{\partial(G\cap\mathbb{C}_I)}ds_I
(s^2-2p_0s+|p|^2)^{-1}(sB-B\overline{p})
\\
&
=
\frac{1}{2\pi}\int_{\partial(G\cap\mathbb{C}_I)}ds_I
(s^2-2p_0s+|p|^2)^{-1}(s-\overline{p})B\\
&
+
\frac{1}{2\pi}\int_{\partial(G\cap\mathbb{C}_I)}ds_I
(s^2-2p_0s+|p|^2)^{-1}(\overline{p}B-B\overline{p})
\end{split}
\]
but observe that
$$
\frac{1}{2\pi}\int_{\partial(G\cap\mathbb{C}_I)}ds_I
(s^2-2p_0s+|p|^2)^{-1}(s-\overline{p})B=
\frac{1}{2\pi}\int_{\partial(G\cap\mathbb{C}_I)}ds_I
S_R^{-1}(s,p)B=B
$$
and moreover by the residue theorem it is
$$
\frac{1}{2\pi}\int_{\partial(G\cap\mathbb{C}_I)}ds_I
(s^2-2p_0s+|p|^2)^{-1}=0
$$
so we get the statement.
\end{proof}

\begin{Tm}\label{PTcommut}
Let $T\in\mathcal{B}^{0,1}(V_n)$ and
 let $\sigma_S(T)= \sigma_{1S}(T)\cup \sigma_{2S}(T)$,
with
$$
{\rm dist}\,( \sigma_{1S}(T),\sigma_{2S}(T))>0.
$$
 Let $U_1$ and
$U_2$ be two axially symmetric s-domains  such that  $\sigma_{1S}(T) \subset U_1$ and $
\sigma_{2S}(T)\subset U_2$,  with $\overline{U}_1
\cap\overline{U}_2=\emptyset$. Set
\begin{equation}\label{pigei}
P_j:=\frac{1}{2\pi }\int_{\partial (U_j\cap \mathbb{C}_I)}S_L^{-1}(s,T) \,
ds_I, \ \ \ \ \ j=1,2,
\end{equation}
\begin{equation}\label{tigei}
T_j:=\frac{1}{2\pi }\int_{\partial (U_j\cap \mathbb{C}_I)}S_L^{-1}(s,T) \,
ds_I\,s, \ \ \ \  j=1,2.
\end{equation}
Then  $P_j$ are projectors and $TP_j=P_jT$ for $j=1,2$.
\end{Tm}
\begin{proof}

Let $\sigma_{jS}(T) \subset G_1$ and $G_2$ be two $T$-admissible open sets such that
$G_1 \cup \partial G_1 \subset G_2 $ and $G_2 \cup \partial G_2 \subset U_j$, for $j=1$ or $2$.
Thanks to the structure of the S-spectrum  we will assume that $G_1$ and $G_2$ are axially symmetric and s-domains.

Take $p\in \partial (G_1\cap \mathbb{C}_I)$ and $s\in \partial (G_2\cap \mathbb{C}_I)$ and observe that, for $I\in \mathbb{S}$,  we have
$$
P_j:=\frac{1}{2\pi }\int_{ \partial ( G_2 \cap \mathbb{C}_I) }  ds_I S^{-1}_R(s,T)
$$
but we can also write $P_j$ as
$$
P_j=\frac{1}{2\pi }\int_{ \partial ( G_1 \cap \mathbb{C}_I) }S^{-1}_L(p,T) dp_I.
$$
Now consider $P_j^2$ written as
$$
P_j^2= \frac{1}{(2\pi)^2 }\int_{ \partial ( G_2 \cap \mathbb{C}_I) }  ds_I \int_{ \partial ( G_1 \cap \mathbb{C}_I) } S^{-1}_R(s,T)S^{-1}_L(p,T) dp_I.
$$
Using the resolvent equation we write:
\[
\begin{split}
P^2&=\frac{1}{(2\pi)^2 }\int_{ \partial ( G_2 \cap \mathbb{C}_I) }  ds_I \int_{ \partial ( G_1 \cap \mathbb{C}_I) }
[S_R^{-1}(s,T)-S_L^{-1}(p,T)]p(p^2-2s_0p+|s|^2)^{-1}dp_I
\\
&
-\frac{1}{(2\pi)^2 }\int_{ \partial ( G_2 \cap \mathbb{C}_I) }  ds_I
\int_{ \partial ( G_1 \cap \mathbb{C}_I) }\overline{s}[S_R^{-1}(s,T)-S_L^{-1}(p,T)](p^2-2s_0p+|s|^2)^{-1}
 dp_I.
 \end{split}
\]
Now observe that
$$
\frac{1}{(2\pi)^2 }\int_{ \partial ( G_2 \cap \mathbb{C}_I) }  ds_I S_R^{-1}(s,T)\int_{ \partial ( G_1 \cap \mathbb{C}_I) }
p(p^2-2s_0p+|s|^2)^{-1}dp_I=0
$$
and
$$
-\frac{1}{(2\pi)^2 }\int_{ \partial ( G_2 \cap \mathbb{C}_I) }  ds_I\overline{s}S_R^{-1}(s,T)
\int_{ \partial ( G_1 \cap \mathbb{C}_I) }(p^2-2s_0p+|s|^2)^{-1}
 dp_I=0
$$
  since the functions
  $$
  p\mapsto p(p^2-2s_0p+|s|^2)^{-1},\ \ \ \ \ \ p\mapsto
  (p^2-2s_0p+|s|^2)^{-1}
  $$ are slice hyperholomorphic and do not have singularities inside $\partial ( G_1 \cap \mathbb{C}_I)$.
 So $P_j^2$ can be written as
\[
\begin{split}
P_j^2&=\frac{1}{(2\pi)^2 }\int_{ \partial ( G_2 \cap \mathbb{C}_I) }  ds_I \int_{ \partial ( G_1 \cap \mathbb{C}_I) }
-S_L^{-1}(p,T)p(p^2-2s_0p+|s|^2)^{-1}dp_I
\\
&
-\frac{1}{(2\pi)^2 }\int_{ \partial ( G_2 \cap \mathbb{C}_I) }  ds_I
\int_{ \partial ( G_1 \cap \mathbb{C}_I) }-\overline{s}S_L^{-1}(p,T)(p^2-2s_0p+|s|^2)^{-1}
 dp_I,\\
& =\frac{1}{(2\pi)^2 }\int_{ \partial ( G_2 \cap \mathbb{C}_I) } \int_{ \partial ( G_1 \cap \mathbb{C}_I) } ds_I
(\overline{s}S_L^{-1}(p,T)-S_L^{-1}(p,T)p)(p^2-2s_0p+|s|^2)^{-1}dp_I.
 \end{split}
\]
Applying now Lemma \ref{Lemma321gf} with  $B:=S^{-1}_L(p,T)$  and observing that $p\in G_2$, we finally have
 $$
 P_j^2=\frac{1}{2\pi }\int_{ \partial ( G_1 \cap \mathbb{C}_I) } S_L^{-1}(p,T) dp_I
=P_j.
$$
Let us now prove that $TP_j=P_jT$. Observe that the functions
$f(s)=s^m$, for $m\in \mathbb{N}_0$ are both right and left
 slice hyperholomorphic. So
 the operator $T$ can be written as
$$
T={{1}\over{2\pi }} \int_{\partial (U\cap \mathbb{C}_I)} S_L^{-1} (s,T)\  ds_I \ s=
{{1}\over{2\pi }} \int_{\partial (U\cap \mathbb{C}_I)} \  s\ ds_I \ S_R^{-1} (s,T);
$$
analogously, as already observed, for the projectors $P_j$ we have
$$
P_j={{1}\over{2\pi }} \int_{\partial (U_j\cap \mathbb{C}_I)} S_L^{-1} (s,T)\
ds_I \ =
{{1}\over{2\pi }} \int_{\partial (U_j\cap \mathbb{C}_I)} \  \ ds_I \ S_R^{-1} (s,T).
$$
From the identity
$$
T_j={{1}\over{2\pi }} \int_{\partial (U_j\cap \mathbb{C}_I)} S_L^{-1} (s,T)\
ds_I \ s={{1}\over{2\pi }} \int_{\partial (U_j\cap \mathbb{C}_I)} \  s\ ds_I \
S_R^{-1} (s,T)
$$
we can compute $TP_j$ as:
$$
TP_j={{1}\over{2\pi }} \int_{\partial (U_j\cap \mathbb{C}_I)} TS_L^{-1} (s,T)\  ds_I \
$$
and using the resolvent equation (\ref{quatSresolrddlftequ})
it follows
\begin{equation}
\begin{split}
TP_j&={{1}\over{2\pi }} \int_{\partial (U_j\cap \mathbb{C}_I)} [S_L^{-1}(s,T)\
s-\mathcal{I}]\  ds_I
\\
&
=
{{1}\over{2\pi }} \int_{\partial (U_j\cap \mathbb{C}_I)} S_L^{-1}(s,T)\ s\  ds_I
\\
&
=
{{1}\over{2\pi }} \int_{\partial (U_j\cap \mathbb{C}_I)} S_L^{-1}(s,T)\  ds_I\ s
\\
&
=T_j.
\end{split}
\end{equation}
Now consider
$$
P_jT={{1}\over{2\pi }} \int_{\partial (U_j\cap \mathbb{C}_I)}
\  \ ds_I \ S_R^{-1} (s,T)T
$$
and using the resolvent equation (\ref{quatSresorigequa})
we obtain
\begin{equation}
\begin{split}
P_jT&={{1}\over{2\pi }} \int_{\partial (U_j\cap \mathbb{C}_I)} \
 \ ds_I \ [s \ S_R^{-1}(s,T)-\mathcal{I}]
\\
&
={{1}\over{2\pi }} \int_{\partial (U_j\cap \mathbb{C}_I)} \
\ ds_I \ s\  S_R^{-1}(s,T)
\\
&
=T_j,
\end{split}
\end{equation}
so the equality $P_jT=TP_j$ holds.
\end{proof}

\begin{Rk}
{\rm
The properties that the Riesz projectors commute with the operator $T$ has been proved for the
 quaternionic version of the S-functional calculus in \cite{acs3}, while the property that $P^2=P$ given in
\cite{css_book} is obtained heuristically. This fact shows the importance of this new resolvent equation.
}
\end{Rk}
As it is well known for hyperholomorphic functions the product of two hyperholomorphic functions is not in general
hyperholomorphic. Here we recall a class of functions for which the pointwise multiplication remains slice hyperholomorphic.
\begin{Dn}\label{defdiN}
Let $f:U\to \mathbb{R}_{n}$ be a slice hyperholomorphic function, where $U$ is an open set in $\mathbb{R}^{n+1}$.
 We define
$$
 \mathcal{N}(U)=\{ f\in\mathcal{SM}(U)\ :  \ f(U\cap \mathbb{C}_I)\subseteq  \mathbb{C}_I,\ \  \forall I\in \mathbb{S}\}.
$$
\end{Dn}

\begin{Pn}\label{monprod}
Let $U$ be an open set in $\mathbb{R}^{n+1}$.
Let $f\in \mathcal{N}(U)$, $g\in \mathcal{SM}(U)$, then $fg\in \mathcal{SM}(U)$.
\end{Pn}

First of all let us observe that  functions in the subclass  $\mathcal{N}(U)$ are both left and right slice hyperholomorphic. When we take the power series expansion of this class of functions at a point on the real line the coefficients of the expansion  are real numbers.

Now observe that for functions in $f\in \mathcal{N}(U)$ we can define $f(T)$
 using the left but also the right S-functional calculus. It is
\[
\begin{split}
 f(T)&={{1}\over{2\pi }} \int_{\partial (U\cap \mathbb{C}_I)} S_L^{-1} (s,T)\  ds_I \ f(s)
 \\
 &
={{1}\over{2\pi }} \int_{\partial (U\cap \mathbb{C}_I)} \  f(s)\ ds_I \ S_R^{-1} (s,T).
\end{split}
\]

\begin{La}\label{Lemma321} Let $B\in \mathcal{B}(V_n)$. Let $G$ be an axially symmetric s-domain and assume that
$f\in \mathcal{N}(G)$.
Then, for $p\in G$, we have
$$
\frac{1}{2\pi}\int_{\partial(G\cap\mathbb{C}_I)}f(s)ds_I
(\overline{s}B-Bp)(p^2-2s_0p+|s|^2)^{-1}=Bf(p).
$$
\end{La}
\begin{proof}
Recalling formula (\ref{efdhsslj}) we write
\[
\begin{split}
\frac{1}{2\pi}\int_{\partial(G\cap\mathbb{C}_I)}& f(s)ds_I
(\overline{s}B-Bp)(p^2-2s_0p+|s|^2)^{-1}
\\
&
=
\frac{1}{2\pi}\int_{\partial(G\cap\mathbb{C}_I)}f(s)ds_I
(s^2-2p_0s+|p|^2)^{-1}(sB-B\overline{p})
\\
&
=
\frac{1}{2\pi}\int_{\partial(G\cap\mathbb{C}_I)}f(s)ds_I
(s^2-2p_0s+|p|^2)^{-1}(s-\overline{p})B\\
&
+
\frac{1}{2\pi}\int_{\partial(G\cap\mathbb{C}_I)}f(s)ds_I
(s^2-2p_0s+|p|^2)^{-1}(\overline{p}B-B\overline{p})
\\
&
:=\mathcal{J}_1+\mathcal{J}_2
\end{split}
\]
but observe that
\[
\begin{split}
\mathcal{J}_1&=\frac{1}{2\pi}\int_{\partial(G\cap\mathbb{C}_I)}f(s)ds_I
(s^2-2p_0s+|p|^2)^{-1}(s-\overline{p})B
\\
&
=
\frac{1}{2\pi}\int_{\partial(G\cap\mathbb{C}_I)}f(s)ds_I
S_R^{-1}(s,p)B=f(p)B.
\end{split}
\]
Consider now the second integral.
 Taking $s=u+Iv$ then the solutions of the equation
 $s^2-2p_0s+|p|^2=0$ are $s_1=\alpha$ and $s_2=\overline{\alpha}$ where $\alpha=p_0+I|\underline{p}|$,
  so
\[
\begin{split}
\mathcal{J}_2&=\frac{1}{2\pi}\int_{\partial(G\cap\mathbb{C}_I)}f(s)ds_I
(s^2-2p_0s+|p|^2)^{-1}(\overline{p}B-B\overline{p})
\\
&
=\frac{1}{2\pi}\int_{\partial(G\cap\mathbb{C}_I)}
\frac{f(s)}{(s-\alpha)(s-\overline{\alpha})}ds_I(\overline{p}B-B\overline{p}),
\end{split}
\]
by the residues theorem we get
\[
\begin{split}
\mathcal{J}_2&=\frac{1}{2\pi}\int_{\partial(G\cap\mathbb{C}_I)}f(s)ds_I
(s^2-2p_0s+|p|^2)^{-1}(\overline{p}B-B\overline{p})
\\
&
=\frac{I}{2|\underline{p}|}[f(p_0-I|\underline{p}|)-f(p_0+I|\underline{p}|)](\overline{p}B-B\overline{p}).
\end{split}
\]
Now we recall the structure formula that shows that
a slice hyperholomorphic function can be written as
$$f(p)=\alpha(p_0,|\underline{p}|)+I_p\beta (p_0,|\underline{p}|)$$
 where
$$
\alpha(p_0,|\underline{p}|)=\frac{1}{2}[f(p_0-I|\underline{p}|)+f(p_0+I|\underline{p}|)],
$$
$$
\beta(p_0,|\underline{p}|)=\frac{I}{2}[f(p_0-I|\underline{p}|)-f(p_0+I|\underline{p}|)]
$$
and in the case of functions $f\in \mathcal{N}$ the functions  $\alpha$ and $\beta$ are real valued. Observe that
\[
\begin{split}
\mathcal{J}_1+\mathcal{J}_2&=f(p)B
+\frac{I}{2|\underline{p}|}[f(p_0-I|\underline{p}|)-f(p_0+I|\underline{p}|)](\overline{p}B-B\overline{p})
\\
&
=\alpha(p_0,|\underline{p}|)B+I_p\beta (p_0,|\underline{p}|)B+\frac{\beta (p_0,|\underline{p}|)}{|\underline{p}|}(\overline{p}B-B\overline{p})
\\
&
=\alpha(p_0,|\underline{p}|)B+I_p\beta (p_0,|\underline{p}|)B+\frac{\beta (p_0,|\underline{p}|)}{|\underline{p}|}((p_0-I_p|\underline{p}|)B-B(p_0-I_p|\underline{p}|))
\\
&
=B(\alpha(p_0,|\underline{p}|)+I_p\beta (p_0,|\underline{p}|))
\\
&
=Bf(p),
\end{split}
\]
so we get the statement.
\end{proof}

 \begin{Rk}{\rm  If we assume that $f\in \mathcal{N}(\mathbf{B}(0,r))$ where
 $\mathbf{B}(0,r)$ is the open ball in $\mathbb{R}^{n+1}$ centered at $0$ and of radius $r>0$  and   $s\in \mathbf{B}(0,r)$, then the proof of the above theorem follows in a shorter way. Indeed we have
$$
(\overline{s}B-Bp)(p^2-2s_0p+|s|^2)^{-1}=\sum_{m\geq 0}s^{-1-m}Bp^m, \ \ \ |p|<|s|.
$$
So
$$
\frac{1}{2\pi}\int_{\partial(G\cap\mathbb{C}_I)}f(s)ds_I\sum_{m\geq 0}s^{-1-m}Bp^m,\ \ \ |p|<|s|,
$$
but
$$
\sum_{m\geq 0}
\frac{1}{2\pi}\int_{\partial(G\cap\mathbb{C}_I)}f(s)ds_Is^{-1-m}Bp^m=
\sum_{m\geq 0}\frac{1}{m!}f^{(m)}(0)Bp^m
$$
and for functions in $\mathcal{N}(\mathbf{B}(0,r))$ the derivatives $f^{(m)}(0)$ are real numbers and so
they commute with $B$. We get
$$
\sum_{m\geq 0}\frac{1}{m!}f^{(m)}(0)Bp^m=B\sum_{m\geq 0}\frac{1}{m!}f^{(m)}(0)p^m=Bf(p).
$$
}
\end{Rk}
We now offer a different proof of the theorem that shows that $(fg)(T)=f(T)g(T)$, under suitable assumptions of $f,g$.
Originally, see \cite{css_book}, the proof was based on the Cauchy formula and the resolvent equations (\ref{quatSresolrddlftequ}),
(\ref{quatSresorigequa}).
\begin{Tm}\label{prodotto}
Let $T\in\mathcal{B}^{\small 0,1}(V_n)$ and assume   $f\in \mathcal{N}_{\sigma_S(T)}$ and $g\in  \mathcal{SM}_{\sigma_S(T)}$.  Then we have
$$
(f g)(T)=f(T)g(T).
$$
\end{Tm}
\begin{proof}
Let $\sigma_{S}(T) \subset G_1$ and $G_2$ be two $T$-admissible open sets such that
$G_1 \cup \partial G_1 \subset G_2 $ and $G_2 \cup \partial G_2 \subset U$.
Take $p\in \partial (G_1\cap \mathbb{C}_I)$ and $s\in \partial (G_2\cap \mathbb{C}_I)$ and observe that, for $I\in \mathbb{S}$,  we have

\[
\begin{split}
 f(T) g(T)&=
 {{1}\over{(2\pi)^2 }} \int_{\partial (G_2\cap \mathbb{C}_I)} \  f(s)\ ds_I \ S_R^{-1} (s,T)
  \int_{\partial (G_1\cap \mathbb{C}_I)} \ S_L^{-1} (p,T)\ dp_I \  g(p)
 \\
&=\frac{1}{(2\pi)^2 }\int_{ \partial ( G_2 \cap \mathbb{C}_I) } f(s)\ ds_I \int_{ \partial ( G_1 \cap \mathbb{C}_I) }
S_R^{-1}(s,T)p(p^2-2s_0p+|s|^2)^{-1}dp_I\  g(p)
 \\
&-\frac{1}{(2\pi)^2 }\int_{ \partial ( G_2 \cap \mathbb{C}_I) } f(s)\ ds_I \int_{ \partial ( G_1 \cap \mathbb{C}_I) }
S_L^{-1}(p,T)p(p^2-2s_0p+|s|^2)^{-1}dp_I\  g(p)
\\
&
-\frac{1}{(2\pi)^2 }\int_{ \partial ( G_2 \cap \mathbb{C}_I) } f(s)\ ds_I
\int_{ \partial ( G_1 \cap \mathbb{C}_I) }\overline{s}S_R^{-1}(s,T)(p^2-2s_0p+|s|^2)^{-1}
 dp_I\  g(p)
 \\
&
+\frac{1}{(2\pi)^2 }\int_{ \partial ( G_2 \cap \mathbb{C}_I) } f(s)\ ds_I
\int_{ \partial ( G_1 \cap \mathbb{C}_I) }\overline{s}S_L^{-1}(p,T)(p^2-2s_0p+|s|^2)^{-1}
 dp_I\  g(p)
 \end{split}
\]
where we have used the resolvent equation. But now observe that
$$
\frac{1}{(2\pi)^2 }\int_{ \partial ( G_2 \cap \mathbb{C}_I) } f(s)\ ds_I \int_{ \partial ( G_1 \cap \mathbb{C}_I) }
S_R^{-1}(s,T)p(p^2-2s_0p+|s|^2)^{-1}dp_I\  g(p)=0
$$
and
$$
-\frac{1}{(2\pi)^2 }\int_{ \partial ( G_2 \cap \mathbb{C}_I) } f(s)\ ds_I
\int_{ \partial ( G_1 \cap \mathbb{C}_I) }\overline{s}S_R^{-1}(s,T)(p^2-2s_0p+|s|^2)^{-1}
 dp_I\  g(p)=0.
$$
so it follows that
\[
\begin{split}
 f(T) g(T)&=
-\frac{1}{(2\pi)^2 }\int_{ \partial ( G_2 \cap \mathbb{C}_I) } f(s)\ ds_I \int_{ \partial ( G_1 \cap \mathbb{C}_I) }
S_L^{-1}(p,T)p(p^2-2s_0p+|s|^2)^{-1}dp_I\  g(p)
\\
&
+\frac{1}{(2\pi)^2 }\int_{ \partial ( G_2 \cap \mathbb{C}_I) } f(s)\ ds_I
\int_{ \partial ( G_1 \cap \mathbb{C}_I) }\overline{s}S_L^{-1}(p,T)(p^2-2s_0p+|s|^2)^{-1}
 dp_I\  g(p)
 \end{split}
\]
which can be written as
\[
\begin{split}
 f(T) g(T)=
 \frac{1}{(2\pi)^2 }\int_{ \partial ( G_2 \cap \mathbb{C}_I) } f(s)\ ds_I &\int_{ \partial ( G_1 \cap \mathbb{C}_I) }
[\overline{s}S_L^{-1}(p,T)-S_L^{-1}(p,T)p]\times
\\
&
\times(p^2-2s_0p+|s|^2)^{-1}dp_I\  g(p).
 \end{split}
\]
Using Lemma \ref{Lemma321} we get
$$
 f(T) g(T)=\frac{1}{2\pi } \int_{ \partial ( G_1 \cap \mathbb{C}_I) }S_L^{-1}(p,T)dp_I \ f(p)\  g(p)
 $$
 which gives the statement.
 \end{proof}
In the original proof of the above theorem we have used the fact that for functions
$f \in \mathcal{N}_{\sigma_S(T)}$ the left S-resolvent equation gives
$$
f(T)T^m=\frac{1}{2\pi } \int_{ \partial ( U \cap \mathbb{C}_I) }S_L^{-1}(p,T)dp_I \ f(p) p^m
$$
from which one obtains
$$
f(T)T^m t^{-1-m}=\frac{1}{2\pi } \int_{ \partial ( U \cap \mathbb{C}_I) }S_L^{-1}(p,T)dp_I \ f(p) p^mt^{-1-m}.
$$
By taking the sum and considering  $t\in \rho_S(T)$, we have
$$
f(T)S_L^{-1}(t,T)=\frac{1}{2\pi } \int_{ \partial ( U \cap \mathbb{C}_I) }S_L^{-1}(p,T)dp_I \ f(p) S_L^{-1}(t,p).
$$
Using this equality and  the Cauchy formula we obtain the statement.

\section{The case of commuting operators and the quaternionic case}

In this last section we state the resolvent equation in the case of commuting operators and for the quaternionic functional calculus.
We also take the occasion to make some comments that show how the
S-functional calculus turns out to be a natural extension of the Riesz-Dunford functional calculus.

\subsection{The case of several commuting operators}

We denote by $\mathcal{BC}^{\small 0,1}(V_n)$ the
subset of  $\mathcal{B}^{\small 0,1}(V_n)$ consisting of paravector operators with commuting components.
 Given an operator in paravector form $T=T_0+e_1T_1+\ldots +e_nT_n$, its so-called conjugate $\overline{T}$ is defined by $\overline{T}=T_0-e_1T_1-\ldots -e_nT_n$. When $T\in \mathcal{BC}^{\small 0,1}(V_n)$ the operator $T\overline{T}$ is well defined and
$T\overline{T}=\overline{T}T=T_0^2+T_1^2+\ldots +T_n^2$ and $T+\overline{T}=2T_0$.

\begin{Tm}\label{lygv}
Let $T\in\mathcal{BC}^{\small 0,1}(V_n)$ and  $s\in \mathbb{R}^{n+1}$ be such that $| s|  < \|T\|$. Then
\begin{equation}\label{SERLEFTQ}
\begin{split}
\sum_{m\geq 0} T^m s^{-1-m}=(s\mathcal{I}-\bar T)(s^2\mathcal{I}-s(T+\overline{T})+T\overline{T})^{-1},
\end{split}
\end{equation}
\begin{equation}\label{SERRIGHQT}
\begin{split}
\sum_{m\geq 0} s^{-1-m} T^m
=
(s^2\mathcal{I}-s(T+\overline{T})+T\overline{T})^{-1}(s\mathcal{I}-\bar T).
\end{split}
\end{equation}
\end{Tm}
The above theorem follows from the fact that the Cauchy kernels for slice hyperholomorphic functions can be written in two possible ways, see Section 2 and \cite{SlaisComm}.
In the case of commuting operators the two expressions are equivalent. The advantage of this approach is that one can work with the so called F-spectrum which  is easier to compute
than the S-spectrum. In fact it can be computed over a complex plane $\mathbb{C}_I$, taking $s=u+Iv$, and then extended to $\mathbb{R}_n$. This is a consequence of the fact that the F-spectrum takes into account the commutativity of the operators $T_j$, $j=0,1,...,n$.
 The F-spectrum is suggested by
 Theorem \ref{lygv} and it is described below.
\begin{Dn}[The $F$-spectrum and the $F$-resolvent sets]\label{FspectrumFres}
Let
$T\in\mathcal{BC}^{\small 0,1}(V_n)$.
We define the $F$-spectrum of $T$  as:
$$
\sigma_F(T)=\{ s\in \mathbb{R}^{n+1}\ \ :\ \ s^2\mathcal{I}-s(T+\overline{T})+T\overline{T}\ \ \
{\it is\ not\  invertible\  }\}.
$$
\noindent
The $F$-resolvent set of $T$ is defined by
$$
\rho_F(T)=\rr^{n+1}\setminus\sigma_F(T).
$$
\end{Dn}
The main properties of the $F$-spectrum are similar to those of the S-spectrum as it is proved in the next results:
\begin{Tm}[Structure of the $F$-spectrum]\label{strutturaK}
 Let $T\in\mathcal{BC}^{\small 0,1}(V_n)$ and let $p = p_0 +p_1I\in
  [p_0 +p_1 I]\subset \mathbb{R}^{n+1}\setminus\mathbb{R}$, such that $p\in \sigma_{F}(T)$.
  Then all the elements of the $(n-1)$-sphere $[p_0 +p_1I]$
 belong to $\sigma_{ F}(T)$.
\end{Tm}

\begin{Tm}[Compactness of $F$-spectrum]
Let
$T\in\mathcal{BC}^{ 0,1}(V_n)$. Then
the $F$-spectrum $\sigma_{F} (T)$  is a compact nonempty set.
Moreover $\sigma_{F} (T)$ is contained in $\{s\in\mathbb{R}^{n+1}\, :\, | s|  \leq \|T\| \ \}$.
\end{Tm}
The relation between the S-spectrum and the F-spectrum is contained in the following result:
\begin{Pn}\label{FeSspettri}
Let $T\in\mathcal{BC}^{\small 0,1}(V_n)$ . Then $\sigma_{F}(T)=\sigma_S(T)$.
\end{Pn}

\begin{Dn}(The $S_C$-resolvent operator)
Let
$T\in\mathcal{BC}^{\small 0,1}(V)$
 and $ s \in \rho_F(T)$.
We define the $S_C$-resolvent operator as
\begin{equation}\label{SCCresolvoperatorQ}
S_{C,L}^{-1}(s,T):=(s\mathcal{I}- \overline{T})(s^2\mathcal{I}-s(T+\overline{T})+T\overline{T})^{-1}.
\end{equation}
\begin{equation}\label{SCCresolvoperatorRQ}
S_{C,R}^{-1}(s,T):=(s^2\mathcal{I}-s(T+\overline{T})+T\overline{T})^{-1}(s\mathcal{I}- \overline{T}).
\end{equation}
\end{Dn}

\begin{Tm}
Let $T\in\mathcal{BC}^{0,1}(V_n)$  and $s, p\in \rho_{F}(T)$. Then
$S_{C,L}^{-1}(s,T)$ satisfies the  left $S$-resolvent equation
\begin{equation}
S_{C, L}^{-1}(s,T)s-TS_{C, L}^{-1}(s,T)=\mathcal{I},
\end{equation}
and $S_{C,R}^{-1}(s,T)$ satisfies
the right $S$-resolvent equation
$$
sS_{C, R}^{-1}(s,T)-S_{C, R}^{-1}(s,T)T=\mathcal{I}.
$$
Moreover, for $p\not\in [s]$, it satisfies the resolvent equation
\[
\begin{split}
S_{C, R}^{-1}(s,T)S_{C, L}^{-1}(p,T)&=((S_{C, R}^{-1}(s,T)-S_{C,L}^{-1}(p,T))p
\\
&
-\overline{s}(S_{C, R}^{-1}(s,T)-S_{C, L}^{-1}(p,T)))(p^2-2s_0p+|s|^2)^{-1}.
\end{split}
\]
which can be written as
\begin{equation}\label{RLresolvIIsc}
\begin{split}
S_{C, R}^{-1}(s,T)S_{C, L}^{-1}(p,T)&=(s^2-2p_0s+|p|^2)^{-1}(s(S_{C, R}^{-1}(s,T)-S_{C, L}^{-1}(p,T))
\\
&
-(S_{C, R}^{-1}(s,T)-S_{C, L}^{-1}(p,T))\overline{p} ).
\end{split}
\end{equation}
\end{Tm}
We conclude this subsection with a couple of considerations on the case of unbounded operators.
\\
(I) Suppose that $T$ is a closed operator with domain $D(T)$.
As one can clearly see, the non commutative version of  $S_L^{-1}(s,T)$, that is
$$
S_L^{-1}(s,T):=-(T^2-2Re(s) T+|s|^2\mathcal{I})^{-1}(T-\overline{s}\mathcal{I}),
$$
is defined on the domain of $T$, and not on $V_n$ as it is in the classical case.
So we have to consider the extension to $V_n$ writing $S_L^{-1}(s,T)$  as follows
$$
\hat{S}_L^{-1}(s,T):=-(T(T^2-2Re(s) T+|s|^2\mathcal{I})^{-1}-(T^2-2Re(s) T+|s|^2\mathcal{I})^{-1}\overline{s}\mathcal{I}).
$$
In this case $\hat{S}_L^{-1}(s,T)$ turns out to be defined on $V_n$.
Observe now that if $T$ is a closed operator with domain $D(T)$ and with commuting components
the left S-resolvent operator $S_{C,L}^{-1}(s,T)$ turns out to be already defined on $V_n$.
For a more detailed discussion see the original papers \cite{formulations, SlaisComm}.
In the case of the right S-resolvent  we have the opposite situation.
With the above consideration the new resolvent equation remains the same also for unbounded operators.
\\
(II) The  F-spectrum is also useful to defined  the so called  F-functional calculus, see \cite{ACSconv, CoSaSo}.
This calculus is defined using the Fueter-Sce-Qian mapping theorem in integral form, see \cite{fueter32, qian, sce}.
It is a hyperholomorphic functional calculus in the spirit of
A.  McIntosh, B. Jefferies and their coauthors (see \cite{jmc,jmcpw,mcp,LMQ}, the monograph \cite{jefferies} and the references therein)
 who first used the theory of hyperholomorphic functions, see \cite{BDS, csss, DSS, ghs}, to define a hyperholomorphic functional calculus for n-tuples of operators.

\subsection{The quaternionic setting}

What we have previously proved in the paper can be rephrased also for the quaternionic functional calculus.
We simply point out that
in this case slice hyperholomorphic functions are defined on an open set $U\subseteq\mathbb{H}$ and have values in the quaternions $\mathbb H$. The resolvent operators are as in the introduction of this paper.
Here it is important to consider right linear operators as well as left linear operators $T$. The possible formulations of the quaternionic functional calculus has been carried out in  \cite{formulations}.
The resolvent equations in Theorem \ref{RLRESOLVEQ} hold in this setting where instead of the paravector operator
$T=T_0+T_1e_1+\ldots +T_ne_n$ we replace quaternionic operators.
We finally mention for sake of completeness one more analogy with the classical case.
As it is well known the Laplace transform of a semigroup $e^{tG}$ where for simplicity we take a bounded operator $G$ defined on a Banach space $X$ is the resolvent operator $(\lambda I-G)^{-1}$. In the quaternionic case we have the analogue result for the two S-resolvent operators:
 Let $T\in \mathcal{B}(V)$ and let $s_0 >\|T\|$.
 Then the left $S$-resolvent operator $S_L^{-1}(s,T)$ is given by
$$
S_L^{-1}(s,T)=\int_0^{+\infty} e^{tT}\, e^{-t s}\, dt,
$$
and $S_R^{-1}(s,T)$ is given by
$$
S_R^{-1}(s,T)=\int_0^{+\infty}  e^{-t s}e^{t\, T} \, dt.
$$
We point out that the theory of the quaternionic evolution operators is developed in \cite{evolution}
where it is also studied the case in which the generator is unbounded.
Recently the case of sectorial operators has been treated  in \cite{GR}.


\begin{thebibliography}{10}



\bibitem{MR2002b:47144}
D. Alpay,
{\em The {S}chur algorithm, reproducing kernel spaces and system
  theory},
 American Mathematical Society, Providence, RI, 2001.
 Translated from the 1998 French original by Stephen S. Wilson,
  Panoramas et Synth\`eses.


\bibitem{acs1}
D.~Alpay, F.~Colombo and I.~Sabadini, {\em Schur functions and their realizations in the slice
hyperholomorphic setting}, Integral
Equations  and Operator Theory {\bf 72} (2012), 253-289.

\bibitem{acs2}
D.~Alpay, F.~Colombo and I.~Sabadini, {\em Pontryagin de
Branges-Rovnyak spaces of slice hyperholomorphic functions}, J.
d'Analyse  Math\'ematique, 2013.

\bibitem{acs3}
D.~Alpay, F.~Colombo and I.~Sabadini, {\em Krein-Langer factorization and related topics in the
slice hyperholomorphic setting},
J. of Geom. Anal., to appear.

\bibitem{ACSconv} D.~Alpay, F.~Colombo and I.~Sabadini, {\em  On some notions of convergence for n-tuples of operators},  Math. Meth. Appl. Sci. (2013/14), to appear.


%
\bibitem{acls}
D.~{Alpay}, F.~{Colombo}, I. Lewkowicz, and I.~{Sabadini}.
\newblock {Realizations of slice hyperholomorphic generalized contractive
and positive functions}.
\newblock arxiv:1310.1035.



\bibitem{adrs}
D.~Alpay, A.~Dijksma, J.~Rovnyak, and H.~de~Snoo.
\newblock {\em {Schur} functions, operator colligations, and reproducing kernel
  {P}ontryagin spaces}, volume~96 of {\em Operator theory: {A}dvances and
  {A}pplications}.
\newblock Birkh{\" a}user Verlag, Basel, 1997.


\bibitem{ALPZ}
 D. Alpay, Y. Peretz, {\rm  Realizations for Schur upper triangular operators centered at an arbitrary point},
  Integral Equations Operator Theory {\bf 37} (2000), 251--323.

\bibitem{av3}
D.~Alpay and V.~Vinnikov,
\newblock  {\rm Finite dimensional de {B}ranges spaces on {R}iemann surfaces}.
\newblock {\em J. Funct. Anal.}, 189(2):283--324, 2002.

\bibitem{BDS} F. Brackx, R. Delanghe, F. Sommen,
\textit{Clifford analysis.} Research Notes in Mathematics, 76, Pitman, Advanced Publishing Program, Boston, MA, 1982.

\bibitem{cs} F. Colombo, I. Sabadini, {\em
On some properties of the quaternionic functional calculus},
J. Geom. Anal., {\bf 19} (2009), 601--627.

\bibitem{evolution} F. Colombo, I. Sabadini,
{\em The quaternionic evolution operator},
 Adv. Math., {\bf 227} (2011), 1772--1805.

\bibitem{formulations} F. Colombo,  I. Sabadini,
{\em On the  formulations of the quaternionic functional calculus},
 J. Geom. Phys.,  {\bf 60} (2010), 1490--1508.


\bibitem{DUKE} F. Colombo,  I. Sabadini,
{\em The Cauchy formula with $s$-monogenic kernel and
 a functional calculus for noncommuting operators},
 J. Math. Anal. Appl., {\bf 373} (2011), 655--679.

\bibitem{SlaisComm}  F. Colombo, I. Sabadini,
{\em The  F-spectrum and the SC-functional calculus},
Proc. Roy. Soc. Edinburgh Sect. A, {\bf  142} (2012),  479--500.

\bibitem{CoSaSo} F. Colombo, I. Sabadini, F. Sommen,
{\em The Fueter mapping theorem in integral form and  the $\mathcal{F}$-functional calculus},
  Math. Meth. Appl. Sci., {\bf 33} (2010), 2050--2066.


\bibitem{functionalcss} F. Colombo, I. Sabadini, D.C. Struppa, {\em
 A new functional calculus for noncommuting operators},
 J. Funct. Anal., {\bf 254} (2008), 2255--2274.


\bibitem{csss} F. Colombo, I. Sabadini, F. Sommen, D.C.
Struppa, {\em Analysis of Dirac Systems and Computational Algebra},
Progress in Mathematical Physics, Vol. 39, {Birkh\"auser}, Boston,
2004.


\bibitem{css_book} F. Colombo, I. Sabadini,  D. C. Struppa,
{\em Noncommutative Functional Calculus. Theory and Applications of Slice Hyperholomorphic Functions},
{Progress in Mathematics}, {\em Birkh\"auser}, Basel,
2011.

\bibitem{DSS} R. Delanghe, F. Sommen, V. Soucek,
{\em Clifford algebra and spinor-valued functions. A function theory for the Dirac operator.}
  Mathematics and its Applications, 53. Kluwer Academic Publishers Group, Dordrecht, 1992. xviii+485 pp.

\bibitem{ds}
N. Dunford, J. Schwartz, {\it Linear Operators, part I: General
Theory }, J. Wiley and Sons (1988).

\bibitem{fueter32}
R. Fueter, {\em Analytische Funktionen einer Quaternionenvariablen}, {Comm. Math.
Helv.}, {\bf 4} (1932), 9--20.

\bibitem{GMP}
R. Ghiloni, V. Moretti, A. Perotti,
{\em Continuous slice functional calculus in quaternionic Hilbert spaces},
Reviews in Mathematical Physics, Vol. 25, No. 4 (2013), pp.1350006--1--1350006--83.

\bibitem{GR} R. Ghiloni, V. Recupero,
{\em Semigroups over real alternative *-algebras: generation theorems and spherical sectorial operators},
Preprint 2013.


\bibitem{ghs} K. G\"urlebeck, K. Habetha, W. Spr\"ossig, { Holomorphic Functions in the Plane and $n$-dimensional space}, Birkh\"auser, Basel, 2008.

\bibitem{jefferies} B. Jefferies, {\em Spectral properties of noncommuting operators},
Lecture Notes in Mathematics, 1843, Springer-Verlag, Berlin, 2004.

\bibitem{jmc} B. Jefferies, A. McIntosh, {\em The Weyl calculus and Clifford analysis},
Bull. Austral. Math. Soc., {\bf 57} (1998), 329--341.

\bibitem{jmcpw} B. Jefferies, A. McIntosh, J. Picton-Warlow, {\em The monogenic functional
calculus}, Studia Math., {\bf 136} (1999), 99--119.

\bibitem{KLMV}
M.S. Liv\u{s}ic, N.~Kravitski, A.~Markus, and V.~Vinnikov.
\newblock {\em Commuting nonselfadjoint operators and their applications to
  system theory}.
\newblock Kluwer, 1995.


\bibitem{mcp} A. McIntosh, A. Pryde, {\em A functional calculus for several commuting
operators}, Indiana Univ. Math. J., {\bf 36} (1987), 421--439.


\bibitem{qian} T. Qian, {\em Generalization of Fueter's result to $\rr^{n+1}$}, Rend. Mat. Acc. Lincei, {\bf 8} (1997), 111--117.
\bibitem{RS} M. Reed, B. Simon, {\em Functional Analysis},
Methods of Modern Mathematical Physics, Academic Press; REV edition, 1980.

\bibitem{LMQ}
C. Li,  A. McIntosh, T. Qian, {\em Clifford algebras, Fourier transforms and singular convolution operators on Lipschitz surfaces},
 Rev. Mat. Iberoamericana {\bf 10} (1994),  665--721.

\bibitem{sce} M. Sce, {\em Osservazioni sulle serie di potenze nei moduli quadratici}, Atti Acc. Lincei Rend. Fisica , {\bf 23} (1957), 220--225.

\bibitem{vinnikov4}
V.~Vinnikov,
\newblock {\em Commuting nonselfadjoint operators and algebraic curves},
  volume~59 of {\em {\rm Operator Theory: Advances and Applications}}, pages
  348--371.
\newblock Birkh{\" a}user Verlag, Basel, 1992.




\end{thebibliography}
\end{document}